\documentclass[11pt]{article}

\usepackage{amsmath, amssymb, amsthm, amstext, graphicx,tikz}
\usepackage{enumerate}
\usepackage{fullpage}
\usetikzlibrary{pgfplots.groupplots}
\tikzstyle{vertex}=[circle,draw=black,fill=black,inner sep=0,minimum size=0.2cm,text=white,font=\footnotesize]
\usepackage[labelfont=bf,labelsep=period]{caption}

\newtheorem{theorem}{Theorem}[section]

\newtheorem{proposition}[theorem]{Proposition}
\newtheorem{definition}[theorem]{Definition}

\newtheorem{conjecture}[theorem]{Conjecture}

\newcommand{\ep}{\varepsilon}

\begin{document}

\title{Robustness of graph properties}
\author{
Benny Sudakov\thanks{Department of Mathematics, ETH, 8092 Zurich, Switzerland.
Email: {\tt benjamin.sudakov@math.ethz.ch}. Research supported by SNSF grant 200021-149111.}}

\date{}

\maketitle

\begin{abstract}
A typical result in graph theory says that a graph $G$, satisfying
certain conditions, has some property $\cal P$. Once such a theorem is
established, it is natural to ask how strongly
$G$ satisfies $\cal P$. Can one strengthen the result by showing that
$G$ possesses $\cal P$ in a robust way? What measures of robustness can one utilize?
In this survey, we discuss various measures that can be used to study robustness of graph
properties, illustrating them with examples.
\end{abstract}

\section{Introduction} 
Let $G$ be a graph and $\mathcal{P}$ a graph property. Many
results in graph theory are of the form ``under certain conditions, $G$
has property $\mathcal{P}$''. Once such a result is established, it
is natural to ask how strongly does $G$ possess $\mathcal{P}$?
In other words, we want to determine the {\em robustness} of $G$ with
respect to $\mathcal{P}$. Recently, there has been increasing interest in the study of robustness of graph properties,
aiming to strengthen classical results in extremal and probabilistic combinatorics.
The goal of this paper is to discuss several such results and to use them to illustrate 
various measures that can be used to study the robustness of graph properties. 

The property which we consider frequently in this survey is Hamiltonicity, which we use as a motivating example.
A \emph{Hamilton cycle} in a graph is a cycle which passes through
every vertex of the graph, and a graph is
\emph{Hamiltonian} if it contains a Hamilton cycle. Hamiltonicity is
one of the most central notions in graph theory which has been intensively
studied by numerous researchers. The problem of deciding the
Hamiltonicity of a graph is one of the NP-complete problems that   
Karp listed in his seminal paper \cite{Karp}, and accordingly, one
cannot hope for a simple classification of such graphs. Nonetheless, there
are many results deriving conditions that are sufficient to establish Hamiltonicity.
For example, a classical result proved by
Dirac in 1952 (see, e.g., \cite[Theorem~10.1.1]{Diestel}) asserts that
every graph on $n \geq 3$ vertices of minimum degree at least $\frac{n}{2}$
is Hamiltonian. In this context, we say that a graph is a {\em Dirac
graph} if it has minimum degree at least $\frac{n}{2}$.
Note that the bound $\frac{n}{2}$ is tight,
as can be seen by the following two examples (in both $n$ is odd): the
first one is a graph obtained by taking two complete
graphs of order $\frac{n+1}{2}$ sharing one vertex,
and the second one is the complete bipartite graph with parts of sizes
$\frac{n+1}{2}$ and $\frac{n-1}{2}$. Both graphs have $n$ vertices and minimum
degree $\frac{n-1}{2}$, but are not Hamiltonian.

Dirac's theorem is one of the most influential results in the study of
Hamiltonicity, and by now many related results are known (see, e.g.,
\cite{Bondy2}). It is therefore very natural to try to strengthen this theorem, by asking
whether Dirac graphs are robustly Hamiltonian. It turns out that there are several ways to answer this question using
different measures of robustness. For example, one could try to show
that a Dirac graph has many Hamilton cycles or that it contains several edge-disjoint Hamilton cycles.
One can also study whether Maker can win a Hamiltonicity game played on the edges of a Dirac graph.
Another natural question concerns the Hamiltonicity of random subgraphs of Dirac graphs. 
We can also put some restrictions on the pairs of edges of a Dirac graph and consider whether there are Hamilton cycles which do not 
contain a pair of conflicting edges. Note that the answer to each of these questions defines in some sense
a measure of robustness of a Dirac graph with respect to Hamiltonicity.

Another measure of how strongly a graph satisfies some property is the so-called {\em resilience} of a property.
A graph property is called {\em monotone increasing} if it is closed under the addition of edges. 
Roughly speaking, for a monotone increasing graph property, the resilience quantifies the
robustness in terms of the number of edges one must delete from $G$,
locally or globally, in order to destroy the property $\mathcal{P}$. 
Resilience was recently studied extensively in the context of random graphs.
The most commonly used model of random graphs, sometimes even synonymous with the term
``random graph", is the so called {\em binomial random graph} $G(n,p)$.

The random graph $G(n,p)$ denotes the probability space whose
elements are graphs on a vertex set $[n]=\{1,\dots, n\}$, and
where each pair of vertices forms an edge, randomly and independently, with
probability $p$. Here, $p$ is a positive real not greater than one, which can depend on $n$.
Abusing notation slightly, we denote a graph drawn from this distribution by $G(n,p)$. It is well known (and easy to prove) that this distribution
is concentrated on the graphs with roughly ${n \choose 2} p$ edges.
We say that the random graph $G(n,p)$
possesses a graph property $\cal P$ {\em asymptotically almost surely},
or a.a.s. for short, if the probability that $G(n,p)$ satisfies $\cal P$
tends to 1 as the number of vertices $n$ tends to infinity. It is well known (see \cite{JLR}) that 
a monotone property $\cal P$  has a {\em threshold} $p_0$ in a sense that if $p \gg p_0$, then  the random graph
a.a.s. satisfies $\cal P$, while if $p \ll p_0$, then $G(n,p)$ a.a.s. does not satisfy it.
The study of random graphs, introduced in the seminal paper of
Erd\H os and R\'enyi \cite{ER60}, has experienced spectacular growth in the last fifty years, with hundreds of 
papers and several monographs (see, e.g., \cite{Bol-book, JLR, FrKar}) devoted to the subject. For almost any interesting graph property, we now understand quite well when the random graph $G(n,p)$ typically has this property. It is therefore interesting to 
obtain robust versions of classical results on random graphs. 

In the subsequent sections of this paper, we discuss the robustness notions  we mentioned above in more detail and 
present several extensions of well-known theorems using these notions. It is of course impossible to cover all 
the known robustness results in one survey, and therefore the choice of
results we present is inevitably subjective. Yet we hope to describe enough examples from this
fascinating area to appeal to many researchers in extremal and probabilistic 
combinatorics, and to motivate further study of the subject.

\section{Many copies}
Given a graph $G$, let $h(G)$ denote the number of distinct Hamilton
cycles in $G$. The above mentioned theorem of Dirac states that
$h(G)>0$ for every $n$-vertex graph with minimum degree at least $n/2$. One obvious way to 
strengthen this theorem is to show that every Dirac graph contains several Hamilton cycles. This was achieved by 
S\'ark\"ozy, Selkow and Szemer\'edi \cite{SSS}, who proved that every Dirac graph $G$ contains
not only one but  at least $c^nn!$ Hamilton cycles for some small
positive constant $c$. They also conjectured that $c$ can be   
improved to $1/2-o(1)$. This is best possible, since  there are Dirac graphs with at most
$(1/2+o(1))^nn!$ Hamilton cycles. To see this, consider the random graph
$G(n,p)$ with $1/2<p=1/2+o(1)$. Indeed, in this case, with high probability the minimum degree
$\delta(G(n,p))=pn+o(n) \geq n/2$ and the expected number of Hamilton cycles
is $p^n(n-1)!$. 

The conjecture of S\'ark\"ozy, Selkow and Szemer\'edi  was proved in a
remarkable paper of Cuckler and
Kahn \cite{CK}. In fact, Cuckler and Kahn proved the following stronger result.
\begin{theorem}
\label{dirac-counting}
For every Dirac graph $G$ on $n$ vertices with minimum degree
$\delta(G)$, 
$$h(G)\geq
\left(\frac{\delta(G)}{e}\right)^{n}(1-o(1))^{n}\,.$$ 
\end{theorem}

\noindent
By considering a random graph with edge probability $p=\delta(G)/n$, one can easily check that this estimate is tight up to a multiplicative factor $(1-o(1))^n$.
The proof of Cuckler and Kahn uses a self-avoiding random walk on $G$, in which the next vertex is chosen 
from the yet unvisited neighbors of the current vertex according to a very cleverly chosen distribution.
The authors show that even after walking for at least $n-o(n)$ steps, one can use the remaining 
vertices to close it into a Hamilton cycle. Since the initial part of the cycle was chosen randomly, this process produces many distinct Hamilton cycles.
The details of how to choose the correct edge-weights and the proof that the above strategy works are quite involved and
will not be discussed here. 

A different approach to this problem was proposed in \cite{FKS}. This approach
is based on the standard estimates for the permanent of a matrix
(the famous Minc conjecture, established by
Br\'egman \cite{Bregman}, and Van der Waerden conjecture, established by Egorychev \cite{Egorychev} and 
Falikman \cite{Falikman}). Let $S_n$ be the set of all permutations of the set $[n]$. 
The \emph{permanent} of an $n\times n$ matrix $A$ is
defined as $per(A)=\sum_{\sigma\in S_n} \prod_{i=1}^n
A_{i\sigma(i)}$. Note that every
permutation $\sigma\in S_n$ has a cycle representation, which is unique up to the order of cycles. 
When $A$ is the  $0$-$1$ adjacency matrix of a graph, every
non-zero summand in the permanent is $1$ and corresponds to a $(\leq 2)$-factor (consisting of cycles and single edges).
Thus,  the permanent of $A$ counts the number of such factors.
A $0$-$1$ matrix $A$ is called $r$-\emph{regular} if it contains exactly $r$ $1$'s in every row and
column. Given such a matrix, the two above mentioned estimates on the permanent show
that 
\begin{eqnarray}
\label{permanent}
r^n \frac{n!}{n^n} \leq  per(A) \leq (r!)^{n/r}, 
\end{eqnarray}
which is asymptotically $(1-o(1))^n (r/e)^n$.

Let $G$ be a graph on $n$ vertices which contains an $r$-factor (that is, an $r$-regular spanning subgraph) with $r$ linear in $n$.
Using estimates (\ref{permanent}) one can show that $G$ contains 
many 2-factors with few (at most $n^{1-\delta}$ for some constant $\delta>0$) cycles. These 2-factors are then converted into
many Hamilton cycles using rotation-extension type  
techniques. For some previous applications of the permanent-based approach to Hamiltonicity problems see, e.g.,  \cite{Alon} and \cite{FK}.

To illustrate this technique we sketch the proof of the following statement,
which gives a lower bound on the number of Hamilton cycles in a
dense graph $G$ in terms of $\textrm{reg}(G)$, where
$\textrm{reg}(G)$ is the maximal even $r$ for which $G$ contains an
$r$-factor. 
\begin{proposition}
\label{weak-dirac-counting}
Let $G$ be a Dirac graph on $n$ vertices. Then the number of Hamilton cycles in $G$ is at
least $\left(\frac{\textrm{reg}(G)}{e}\right)^n (1-o(1))^n$.
\end{proposition}

\noindent
In particular, for a $cn$-regular graph $G$ on $n$ vertices and constant $c>1/2$
we show that $h(G) \geq (c+o(1))^n n!$, which is  asymptotically tight.
In general, the estimate in Proposition \ref{weak-dirac-counting} is weaker than the result of Cuckler and Kahn.
On the other hand, since every Dirac graph contains an
$r$-factor with even $r$ about $n/4$ (see \cite{KAT}), this bound implies immediately
the result of S\'ark\"ozy, Selkow and Szemer\'edi mentioned above. 

\vspace{0.2cm}
\noindent
{\bf Sketch of proof of Proposition \ref{weak-dirac-counting}.}\, 
Let $r=\textrm{reg}(G)$ and let $H\subseteq G$ be an 
$r$-factor of $G$. We have that $r=\Theta(n)$ and we can use permanent estimates (see \cite{FKS}) to show that
$H$ has at least $(r/e)^n (1-o(1))^n$ $2$-factors 
with at most $s=n^{1/2+o(1)}$ cycles. For every such 
factor $F$, we can
turn $F$ into a Hamilton cycle of $G$ by adding and removing at most
$O(s)$ edges as follows.

Take a non-Hamilton cycle $C$ in $F$. By the connectivity of $G$ we can find a vertex $v\in V(C)$ and a vertex $u\in
V(G)\setminus V(C)$ for which $vu\in E(G)$. Then by deleting the
edge $vv^+$ from $C$ ($v^+$ is next vertex after $v$ in the cycle) we get a path $P$ which can be extended by the edge   
$vu$. Connecting it to a cycle $C'$ which contains $u$ we obtain a longer path $P'$. Repeat this argument as
long as we can. If there are no edges between the endpoints $w, w'$ of the
current path $P'$ and the other cycles from $F$, then we can close $P'$ into a cycle as in the standard proof of Dirac's theorem.
Indeed, both $w, w'$ have at least $n/2$ neighbors in $P'$ and therefore there is a pair of consecutive vertices $v, v^+$ of $P'$
such that $w$ is adjacent to $v^+$ and $w'$ is adjacent to $v$. 

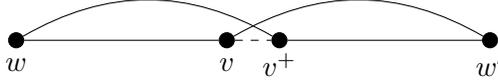
\begin{figure}[h]
\begin{center}
\begin{tikzpicture}[scale=1.4]

\def\labelsep{0.23}
\node[vertex] (u) at (0,0) {};
\node[vertex] (vi) at (2,0) {};
\node[vertex] (vi1) at (2.5,0) {};
\node[vertex] (v) at (4.5,0) {};
\node at (0,-\labelsep) {$w$};
\node at (2,-\labelsep) {$v$};
\node at (2.5,-\labelsep) {$v^+$};
\node at (4.5,-\labelsep) {$w'$};

\draw (u) -- (vi) to [bend left] (v);
\draw (v) -- (vi1);
\draw (vi1) to [bend right] (u);
\draw [dashed] (vi) -- (vi1);

\end{tikzpicture}
\end{center}
\end{figure}

\vspace{-0.7cm}
\noindent
These edges together with $P'$ form a cycle,
which is connected to some other cycle in $F$ as we explained above. 
Note that in each such step we invest at most $4$ edge replacements
in order to decrease the number of cycles in the factor by $1$, and unless the  
current cycle is a Hamilton cycle, we can always merge two cycles.
Therefore, after $O(s)$ edge replacements we get a Hamilton cycle.

In order to complete the proof, note that every Hamilton cycle $C$ we constructed was counted 
at most $\binom{n}{s}(2s)^{2s}=r^{o(n)}$ times. Indeed, to get a 2-factor $F$ from $C$, choose $s$ edges of $C$
to delete. This gives at most $s$ paths which need to be turned into a
$2$-factor by connecting their endpoints. For each endpoint we
have at most $2s$ choices of other endpoints to connect it to.
Therefore, we obtain at least $(r/e)^n (1-o(1))^n$ Hamilton cycles.
\hfill $\Box$

\vspace{0.2cm}

A {\em tournament} is an oriented complete graph. It is easy to prove by induction that every 
tournament has a Hamilton path. Moreover, if the tournament is regular, i.e., if all vertices have the same in/outdegrees, then
it also has a Hamilton cycle. Hence, it is logical to ask whether it has many such cycles.
Counting Hamilton cycles in tournaments is a very old problem
which goes back some seventy years to one of the first applications
of the probabilistic method by Szele \cite{Sz}. He proved
that there is a tournament on $n$ vertices which has at least $(n-1)!/2^n$ Hamilton
cycles. Alon \cite{Alon} showed, using permanent estimates, that this result is nearly tight and that
every $n$-vertex tournament has at most $O(n^{3/2}(n-1)!/2^n)$   
Hamilton cycles.  Thomassen \cite{Tho1} conjectured that the randomness is unnecessary in Szele's
result and that in fact every regular tournament contains at least
$n^{(1-o(1))n}$ Hamilton cycles. This conjecture was solved by
Cuckler \cite{Cuckler} using the above mentioned random walk approach. He 
proved that every regular tournament on $n$ vertices contains at least $\frac{n!}{(2+o(1))^n}$ Hamilton
cycles. 
 
Tournaments are a special case of \emph{oriented} graphs, i.e., directed graphs which are obtained by orienting the 
edges of a simple graph. That is, between every unordered pair of
vertices $\{x,y\}\subseteq V(G)$ at most one of the
(oriented) edges $xy$ or $yx$ is present. Hamiltonicity problems in oriented 
graphs are usually much more challenging. Given an oriented graph  
$G$, let $\delta^+(G)$ and $\delta^-(G)$ denote the minimum
\emph{outdegree} and \emph{indegree} of the vertices in $G$,
respectively. We also use the notation $\delta^{\pm}(G)=\min\{\delta^+(G),\delta^-(G)\}$ and refer to it as
the \emph{minimum semi-degree} of $G$. In the late 70's Thomassen \cite{Tho}
raised the question of determining the minimum semi-degree
that ensures the existence of a Hamilton cycle in an oriented graph 
$G$. H\"aggkvist \cite{Haggkvist} found a construction which gives a
lower bound of $\frac{3n-4}{8}-1$. The problem was only resolved 
recently by Keevash, K\"uhn and Osthus \cite{KKO}, who proved an analogue of Dirac's theorem for oriented and large enough graphs.
\begin{theorem}
\label{oriented-dirac}
For all sufficiently large $n$, every oriented graph $G$ on $n$ vertices with $\delta^{\pm}(G)\geq
\frac{3n-4}{8}$ contains a Hamilton cycle.
\end{theorem}
Having obtained such a theorem, it is in the spirit of the above discussion to ask whether one can actually find many 
Hamilton cycles when $\delta^{\pm}(G)\geq \frac{3n-4}{8}$. It turns out that the permanent-based approach 
can be used to tackle this question in the case the oriented graph is nearly regular, yielding the following theorem \cite{FKS}.

\begin{theorem} \label{CountingHamOriented}
Let $n$ be sufficiently large and let $G$ be an oriented graph on $n$ vertices whose 
in/outdegrees are all $cn\pm o(n)$ for some constant $c>3/8$.
Then $h(G)\geq  \left(\frac{(c+o(1))n}{e}\right)^n$. 
\end{theorem}  

\noindent
The bound on the in/outdegrees in this theorem is tight. This follows from
the above-mentioned construction of  H\"aggkvist \cite{Haggkvist}, which shows that there are non-Hamiltonian $n$-vertex oriented graphs with
all in/outdegrees of order $(3/8-o(1))n$. Since in a regular tournament 
all in/outdegrees are $\frac{n-1}{2}$, this theorem substantially extends the result of Cuckler \cite{Cuckler} mentioned above.

The Hamiltonicity problem was extensively studied for random  graphs as well.
Early results on the Hamiltonicity of random graphs were proved by P\'osa \cite{Posa}
and Korshunov \cite{Korshunov}. Improving on these
results, Bollob\'as \cite{Bollobas84} and Koml\'os
and Szemer\'edi \cite{KoSz} proved that if $p\ge(\log n+\log\log n+\omega(n))/n$
for any function $\omega(n)$ that goes to infinity together with  
$n$, then $G(n,p)$ is a.a.s. Hamiltonian (here and later in the paper all logarithms are natural). The range of $p$ cannot
be improved, since if $p\le(\log n+\log\log n-\omega(n))/n$, then
$G(n,p)$ a.a.s. has a vertex of degree at most one.  

Once the Hamiltonicity threshold is established, it is natural to 
try to estimate the number of Hamilton cycles in the random graph.
Using linearity of expectation, one can immediately see that 
the expected number of such cycles in $G(n,p)$ is $\frac{(n-1)!}{2}p^n$. 
Therefore it is logical to suspect that the actual number 
of Hamilton cycles in the random graph is typically close 
to this value, at least when the edge probability is not too small.
This was indeed confirmed by Glebov and Krivelevich \cite{GK} (also see their paper for the history of the problem and many related results).
They proved that for all edge probabilities for which $G(n,p)$ is Hamiltonian, i.e., for $p\ge(\log n+\log\log n+\omega(n))/n$, it has 
a.a.s. $n!p^n (1-o(1))^n$ Hamilton cycles.

\section{Edge disjoint copies}
In the previous section we showed that any Dirac graph contains not only one but exponentially many Hamilton cycles. 
Another interesting way to extend Dirac's theorem is to show that a minimum degree of at least $n/2$ implies the existence of many edge-disjoint Hamilton cycles. Note that there is no obvious way to deduce such a statement from the
results in Section 2. Indeed, a graph might have many Hamilton cycles, all sharing the same small set of edges. The problem of 
how many edge-disjoint Hamilton cycles one can find in a Dirac graph was first posed by 
Nash-Williams \cite{NashWilliams1} in 1970. In \cite{NashWilliams2} he proved that any such graph has
 at least $\frac{5}{224}n$ edge-disjoint Hamilton cycles. He also asked \cite{NashWilliams1, NashWilliams2} to improve this estimate. Clearly,
 $\lfloor(n+1)/4\rfloor$ is a general upper bound on the number of edge-disjoint Hamilton cycles in a Dirac graph obtained by considering an $n/2$-regular graph, and originally
 Nash-Williams \cite{NashWilliams1} believed this to be tight.
 
 Babai (see also \cite{NashWilliams1}) found a counterexample to this conjecture. Extending his ideas further, Nash-Williams gave an example of a graph on $n=8k+2$ vertices with
 minimum degree $4k+1$ and with at most $(n-2)/8$ edge-disjoint Hamilton cycles. He conjectured that this example is tight, that is, any Dirac graph contains at least
 $(n-2)/8$ edge-disjoint Hamilton cycles. Recall that $\textrm{reg}(G)$ is the maximal even $r$ for which $G$ contains an
 $r$-factor. Since a Hamilton cycle is a 2-regular subgraph of $G$, the maximum number of edge-disjoint Hamilton cycles in $G$ is at most
 $\textrm{reg}(G)/2$. In order to study the maximum number of edge-disjoint Hamilton cycles in Dirac graphs, it is therefore natural to ask for the largest even integer $r$ such 
 that every $n$-vertex graph with minimum degree $\delta$ must contain an $r$-regular spanning subgraph. We denote this function by $\textrm{reg}(n,\delta)$. 
 Note that the complete bipartite graph whose parts differ by one vertex shows that $\textrm{reg}(n,\delta)=0$ for $\delta<n/2$.
 The case $\delta=n/2$ was solved by Katerinis \cite{KAT}, who showed that the above mentioned example of 
 Nash-Williams is tight and $\textrm{reg}(n,n/2) = (n-2)/8$. For $\delta>n/2$ it is known (see, e.g., \cite{CKLOT}) that 
 $$\frac{\delta+\sqrt{n(2\delta-n)}}{2}-1 \leq \textrm{reg}(n,\delta) \leq \frac{\delta+\sqrt{n(2\delta-n)}}{2} + 1.$$
 
 After many partial results, the question of Nash-Williams was answered asymptotically by Csaba, K\"uhn, Lo, Osthus and Treglown \cite{CKLOT}. They proved that for large $n$, every $n$-vertex Dirac graph $G$ contains at least
 $\textrm{reg}(n,\delta(G))/2$ edge-disjoint Hamilton cycles. Although this determines the worst-case behavior,  for a specific Dirac graph $G$, $\textrm{reg}(G)$ can be much larger than $\textrm{reg}(n,\delta)$. Therefore, it is logical to try and bound the maximum number of edge-disjoint Hamilton cycles in terms of $\textrm{reg}(G)$. This question was raised by 
 K\"uhn, Lapinskas and Osthus \cite{KLO}, who conjectured the following tight result.
 
 \begin{conjecture}\label{regconj}
 	Suppose $G$ is a Dirac graph. Then $G$ contains at least
 	$\textrm{reg}(G)/2$ edge-disjoint Hamilton cycles.
 \end{conjecture}

 The following theorem of Ferber, Krivelevich and the author \cite{FKS} gives an approximate asymptotic version of this conjecture.
 
 \begin{theorem}\label{AppRegConj} For every $\varepsilon>0$ and a sufficiently large integer $n$, the following holds: every graph $G$ on $n$ vertices and with $\delta(G)\geq (1/2+
 \varepsilon)n$ contains at least $(1-\varepsilon)\textrm{reg}(G)/2$ edge-disjoint Hamilton cycles. \end{theorem}
 
 \noindent
{\bf Sketch of proof.}\, The approach based on the permanent estimates and rotations, described in the previous section, is rather oblivious to the value of
 $\textrm{reg}(G)$. This makes it well-suited for the proof of the above theorem.
 
We first construct an auxiliary graph $H \subset G$ with relatively small degree and good expansion properties which we use to perform rotations.
More precisely, we claim that for every constant $\varepsilon$ and $\alpha \leq \varepsilon^2$ there exist some $\beta \ll \varepsilon \alpha$ and a subgraph  $H \subset G$ satisfying all of the following: $\delta(H) \geq \varepsilon n/8$ and $G-H$ still has an $r$-factor $G'$ for an even $r$ with $r \geq (1-\varepsilon/2)\textrm{reg}(G)$. Moreover,
for every subset of edges $E'$ such that $|E'| \leq \beta n^2$ and  $\delta(H-E') \geq \alpha n$ we have that $H-E'$ is connected and every subset $S, |S| \geq \alpha n$ in $H-E'$ 
has at least
$(1/2+\varepsilon/4)n$ neighbors outside $S$. This graph can be obtained by randomly choosing $\varepsilon n/16$ 2-factors in $G$. To prove that it has all the listed above properties 
we use the permanent estimates from Section 2. Let $S,T\subseteq V(G)$ be two disjoint subsets of sizes $|S|=\alpha n$ and $|T|=\frac{(1-\varepsilon)n}{2}$. 
It is enough to show that with high probability $|E_{H}(S,T)|\geq \beta n^2$.
Since $\delta(G)\geq (1/2+\varepsilon)n$, it follows that $d_G(v,T)\geq \varepsilon n/2$ for every $v\in S$. Therefore, $|E_G(S,T)|\geq
|S|\varepsilon n/2=\frac{\varepsilon\cdot \alpha}{2} n^2 \gg \beta n^2$. Using this fact one can bound the probability that 
$|E_{H}(S,T)|\geq \beta n^2$ as a ratio of $per(A')/per(A)$,  where $A$ is the adjacency matrix $G$ and $A'$ is obtained from $A$ by deleting
at least $\frac{\varepsilon\cdot \alpha}{2} n^2 -\beta n^2 \geq \frac{\varepsilon\cdot \alpha}{4} n^2$ ones. The permanent estimates can be used to show that this ratio is exponentially small so that the union bound over all
pairs $S, T$ can be applied (see \cite{FKS} for more details).

Now the proof can be completed similarly to that of Proposition \ref{weak-dirac-counting}.
We repeatedly find and delete  $m=(1-\varepsilon/2)r/2$
edge-disjoint $2$-factors of $G'$, each containing at most $s^*=n^{1/2+o(1)}$ cycles. 
Note that by removing such a factor from an
$r'$-regular graph, the resulting graph is $(r'-2)$-regular, and
therefore one can apply the permanent estimates over and over again. 
Finally we perform rotations using the edges of $H$ to turn every such factor into a Hamilton cycle
by replacing at most $O(s^*)$ edges. Moreover, every edge of $H$ which we use is  permanently deleted from $H$.
Since the total number of deleted edges in all steps is at most 
$O(n s^*)=o(n^2)$, the properties of $H$ are not affected.
\hfill $\Box$

\vspace{0.2cm}
The example of a Dirac graph with at most $(n-2)/8$ edge-disjoint Hamilton cycles
is not regular. Accordingly, Nash-Williams conjectured that
every $d$-regular Dirac graph contains $\lfloor d/2\rfloor$ edge-disjoint Hamilton cycles. Recently, in a remarkable tour de force,  Csaba, K\"uhn, Lo, Osthus and Treglown \cite{CKLOT} proved this conjecture for all sufficiently large Dirac graphs.

Recall that a tournament is an orientation of the complete graph. In the previous section we mentioned that every  regular tournament has 
a Hamilton cycle. Does it have many edge-disjoint Hamilton cycles? In principle, since all in/outdegrees are equal it is possible that it contains $(n-1)/2$ edge-disjoint Hamilton cycles. This question was raised by Kelly in 1968, who posed the following striking conjecture.

\begin{conjecture}
The edges of every regular tournament  can be decomposed into Hamilton cycles.
\end{conjecture}

Recently, this was proved for large tournaments by K\"uhn and Osthus \cite{KO1}. In order to solve the conjecture they developed a very powerful decomposition theorem 
for the so called {\em robust outexpanders}. Given an $n$-vertex directed graph $G$ a {\em $\nu$-robust outneighborhood} of a vertex-set $S$ is the set of all vertices in $G$ with at least $\nu n$ inneighbors in $S$.
The graph is called a {\em robust $(\nu,\tau)$-outexpander} if the $\nu$-robust outneighborhood of every set $S$ of size $\tau n \leq  |S| \leq (1-\tau)n$ has size at least $|S|+\nu n$.
We call an oriented graph $r$-regular if all its in/outdegrees are equal to $r$. Using the celebrated Szemer\'edi's Regularity Lemma (see, e.g., \cite{KomSi} for many other applications of this important tool), 
K\"uhn and Osthus \cite{KO1} proved that a large regular robust outexpander with linear degree has a Hamilton decomposition. 
One consequence of their result is the following decomposition theorem for regular oriented graphs of high degree, see \cite{KO2}. 
\begin{theorem}
For all sufficiently large $n$ and constant $c>3/8$, every $cn$-regular oriented $n$-vertex graph $G$ has a Hamilton decomposition.
\end{theorem}

\noindent
To prove this theorem, it is enough to simply verify that such oriented graphs are robust outexpanders. The constant $3/8$ is tight, since as we mentioned in Section 2, there are $cn$-regular oriented $n$-vertex graphs with $c<3/8$ which are not Hamiltonian. Since a regular tournament is an $(n-1)/2$-regular oriented graph, this theorem implies Kelly's conjecture
for large tournaments. 

An oriented $n$-vertex graph with minimum semi-degree at least $3n/8$ is Hamiltonian by Theorem \ref{oriented-dirac}.
Can one strengthen this result by showing that such a graph has many edge-disjoint Hamilton cycles? Similarly to the undirected case,  for an oriented graph 
$G$ let $\textrm{reg}(G)$ be the maximum integer $r$ such that $G$ has a spanning $r$-regular subgraph. If $G$ has $t$ edge-disjoint Hamilton cycles then clearly their union is $t$-regular and therefore $t \leq \textrm{reg}(G)$.
Together with Ferber and Long \cite{FLS}, we believe that this bound is tight.
\begin{conjecture}
Let $n$ be sufficiently large, and let $G$ be an oriented graph on $n$ vertices with $\delta^\pm(G)\geq 3n/8$. Then $G$ contains $\textrm{reg}(G)$ edge-disjoint Hamilton cycles.
\end{conjecture}

\noindent
In \cite{FLS} we prove the following approximate version of this conjecture. If $G$ is an $n$-vertex oriented graph with $\delta^\pm(G)\geq  cn$ for some constant $c>3/8$ and $n$ is sufficiently large, then  $G$ contains $(1-o(1))reg(G)$ edge-disjoint Hamilton cycles. 
The proof is based on various probabilistic arguments and gives a rather short alternative proof of an approximate version of Kelly's conjecture, first established in \cite{KOT}.

Given a regular graph $G$, denote by $H(G)$ the number of Hamiltonian decompositions of $G$. 
We already know that such a $G$ of high degree has a Hamiltonian decomposition, so it is interesting to count the number of decompositions.
To estimate the number of Hamilton cycles in $G$, observe that $per(A_G)$ is an upper bound. Combining this with
the permanent estimates explained in (\ref{permanent}), we see that an $r$-regular graph has at most
$(r!)^{n/r}$ Hamilton cycles. By choosing any Hamilton cycle and  deleting it, we obtain an $(r-2)$-regular graph. The number of Hamilton cycles in the new graph can again be bounded from above using (\ref{permanent}).
Continuing this process and multiplying all the estimates for regularities $r, r-2, r-4, \ldots $
we can use Stirling's formula to deduce that $H(G)$ is at most 
$\left((1+o(1))\frac{r}{e^2}\right)^{nr/2}$. In the case of an oriented $r$-regular graphs the same arguments gives that the number of Hamilton decompositions is at most $\left((1+o(1))\frac{r}{e^2}\right)^{nr}$. 
It turns out that some of the tools developed to prove the existence of many edge-disjoint Hamilton cycles can be pushed further to count the decompositions.
In particular, for sufficiently dense graphs, both of these estimates have the right order of magnitude. This was established for $r$-regular $n$-vertex graphs with $r \geq (1/2+\epsilon)n$ in \cite{GLS}  and for 
oriented graphs with $r \geq (3/8+\epsilon)n$ in \cite{FLS}.

As we explained above, the problem of determining the number of edge-disjoint Hamilton cycles that can be packed into a given graph has a long history and was extensively studied for various classes of graphs.
We conclude this section by discussing this question in the context of random graphs. In order to contain $s$ edge-disjoint Hamilton cycles the graph must clearly have minimum degree of at least $2s$. For $G(n,p)$, this happens typically when
$p\ge(\log n+(2s-1)\log\log n+\omega(n))/n$. Generalizing such Hamiltonicity result for random graphs, 
Bollob\'as and Frieze \cite{BF} proved that a minimum degree of $2s$ is indeed a.a.s. sufficient for $G(n,p)$ to contain $s$ edge-disjoint Hamilton cycles. This motivated Frieze and Krivelevich \cite{FK2} to conjecture that a random graph $G=G(n,p)$ a.a.s. has $\lfloor \delta(G)/2 \rfloor$ edge-disjoint Hamilton cycles for all edge probabilities. This striking conjecture was confirmed in a series of papers by several researches, with two main ranges of edge probabilities 
covered in \cite{KnKO, KS}. The problems of packing and counting Hamilton cycles 
in random directed graphs were studied in \cite{FKL}.

\section{Random subgraphs and Maker-Breaker games}
In this section we discuss two additional, closely related, robust extensions of Dirac's theorem.
In both cases, instead of having access to all the edges of a Dirac graph $G$, we have to find a 
Hamilton cycle in some very sparse subgraph of $G$. Moreover, this subgraph is given to us by some random process 
or by an adversary.

\subsection{Random subgraphs}
An equivalent way of describing the random graph $G(n,p)$ is to say that it is the probability space
of graphs obtained by taking every edge of the complete graph $K_n$ 
independently with probability $p$. A variety of questions can be  
asked when we start with a host graph $G$ other than $K_n$, and
consider the probability space of graphs obtained by taking each one of its  
edge independently with probability $p$. We denote this  
probability space as $G_p$. If the original graph $G$ has property $\cal P$, then a natural way to 
strengthen such a statement is to show that for some $p<1$ the random subgraph $G_p$ with high probability 
also has this property. By doing so, one frequently obtains interesting and
challenging questions.

For example, it is an easy exercise to prove that a graph $G$ with minimum degree $k$ contains a path and even a cycle 
with at least $k+1$ vertices. On the other hand, it is more difficult to show that this remains true for a random subgraph
$G_p$ of such $G$ when $p<1$. An even more challenging question is to determine all values of $p$ such that $G_p$ with high probability has a cycle of length at least $k+1$.
This and the related questions were recently studied by several researchers  in \cite{ KLS1, Rior, GNS}. In \cite{GNS} it was  proved that for large $k$ and 
$p\geq \frac{\log k+\log \log k+\omega(k)}{k}$ the random graph $G_p$ a.a.s. contains a cycle of length at least $k+1$. Since one can take $G$ to be a complete graph on $k+1$ vertices, this result
generalizes the classical result on Hamiltonicity of the random graph $G(k,p)$.

In this context, to get a better understanding of the robustness of Dirac's
theorem we consider the following question.  Let $G$ be an $n$-vertex graph with minimum degree at least $\frac{n}{2}$. Since Hamiltonicity is a monotone graph property, we know (see \cite{JLR}) that there
exists a threshold $p_0$ such that if $p \gg p_0$,
then $G_p$ is a.a.s. Hamiltonian, and if $p \ll p_0$, then it is
a.a.s. not Hamiltonian. What is the Hamiltonicity threshold for $G_p$, in
particular, does $G_p$ stay Hamiltonian for $p \ll 1$? Note that a positive 
answer to the latter question shows  that typically one cannot destroy Hamiltonicity even 
by randomly removing most of the edges of a Dirac graph. The following theorem, which was proved in \cite{KLS2}, 
answers these questions.

\begin{theorem} \label{Dirac-random-subgraph}
There exists a positive constant $C$ such that for $p\ge\frac{C\log
n}{n}$ and a Dirac graph $G$ on $n$ vertices, the random subgraph $G_{p}$ is a.a.s. Hamiltonian.
\end{theorem}

\noindent
This theorem establishes the correct order of magnitude of the  
threshold function  for Hamiltonicity of the random subgraph of any Dirac graph. Indeed, if $p \le (1+o(1))\frac{\log n}{n}$, then
it is easy to see that the graph $G_p$ a.a.s. has isolated vertices. 

It is worth comparing this theorem with the robustness
results from the two previous sections. Given a Dirac graph $G$, recall that $h(G)$ is the number
of Hamilton cycles in $G$. Since the expected number of Hamilton cycles in the random subgraph 
$G_p$ is $p^n h(G)$, Theorem \ref{Dirac-random-subgraph} implies that $p^n h(G) \ge 1$
for $p \ge \frac{C\log n}{n}$. This shows that $h(G) \ge \left(\frac{n}{C\log n}\right)^n$.
We can also partition the edges of a Dirac graph $G$ into $t=C^{-1} n/\log n$ disjoint parts, by putting every edge 
randomly and independently into one of these parts with probability
$p=1/t$. Every part than behaves like a random subgraph $G_p$ and a.a.s. contains a Hamilton cycle. We have therefore obtained a partition such that most of its parts contain a 
Hamilton cycle, showing that a Dirac graph has at least $\Omega(n/\log n)$ edge-disjoint Hamilton cycles.
The above theorem can thus  be used to recover somewhat weaker versions of the results mentioned in Sections 2 and 3.

\vspace{0.2cm}
\noindent
{\bf Sketch of proof of Theorem \ref{Dirac-random-subgraph}.}\, 
The proof is rather involved, so we only give a very brief overview.
Our main tool is P\'osa's rotation-extension technique (see \cite{Posa}),
which exploits the expansion property of the graph. Let $G$ be a connected graph and let $P = (v_0, \cdots, v_\ell)$ be a path in $G$ which we want to extend. Suppose that we cannot do this directly, i.e., all neighbors of $v_0, v_\ell$ lie on $P$.
If $G$ contains an edge of the form $(v_0, v_{i+1})$ for some $i$,
then $P' = (v_i,\cdots, v_0, v_{i+1},\cdots, v_{\ell})$
forms another path of length $\ell$ in $G$. 
\begin{figure}[h]
  \centering\begin{center}
  \begin{tikzpicture}[scale=1.5]
  
  \def\labelsep{0.23}
  \node[vertex] (u) at (0,0) {};
  \node[vertex] (xm) at (2,0) {};
  \node[vertex] (x) at (2.5,0) {};
  \node[vertex] (v) at (4.5,0) {};
  \node at (0,-\labelsep) {$v_0$};
  \node at (2,-\labelsep) {$v_i$};
  \node at (2.5,-\labelsep) {$v_{i+1}$};
  \node at (4.5,-\labelsep) {$v_\ell$};
  \node at (3.5,0.8*\labelsep) {$P$};
  
  \draw (xm) -- (u) to [bend left] (x);
  \draw (x)--(v);
  \draw [dashed] (xm) -- (x);
  
  \end{tikzpicture}
  \end{center}
  \end{figure}
We say that $P'$ is obtained from $P$ by a
\emph{rotation} with {\em fixed endpoint} $v_\ell$, \emph{pivot point}
$v_{i+1}$ and {\em broken edge} $(v_i, v_{i+1})$. Note that after
performing this rotation, we can now try to extend the new path using edges incident to $v_i$.
We can also close the cycle using the edge $(v_i, v_\ell)$, if it exists.
Since $G$ is connected this will also allow us to extend the path.
As we perform more and more rotations, we will get more and more new endpoints and such
closing pairs. We employ this rotation-extension technique repeatedly until we can extend the path.

The proof of Theorem \ref{Dirac-random-subgraph} splits into three different cases, based on the structure of the Dirac graph $G$.
The first two cases apply when $G$ is close to one of the two extremal configurations, i.e., if it is close to a complete bipartite graph or to the union of two disjoint
complete graphs. These cases are easier to handle since one can use techniques developed to show Hamiltonicity of random graphs.
The third case applies when $G$ has at least $cn^2$ edges between any two disjoint subsets of size $n/2$ for some constant $c>0$.
In this case one can prove that a random subgraph $G_p$ has typically strong expansion properties and in particular still has many edges
between any two disjoint subsets of size $n/2$. Using this expansion property of $G_p$ one can carefully perform rotations of its longest path
to obtain a set $S_P$ of linear size with the following property: for every $v \in S_P$ there is a set $T_v$ of size at least $(1/2+\delta)n$ with constant $\delta>0$ such that for every
$w \in T_v$ there is a path of maximal length from $v$ to $w$. Since the degree of $v$ in $G$ was at least $n/2$, there are linearly many edges in $G$ between $v$ and $T_v$. 
If any of these edges is in $G_p$ we can close the longest path into a cycle, which is a Hamilton cycle by maximality of the initial path and the connectivity of $G_p$ . Since the set $S_P$ has linear size, one can easily show that this a.a.s. occurs for at least one vertex in $S_P$. \hfill $\Box$

\subsection{Maker-Breaker games} 
Let $V$ be a set of elements and $\mathcal{F} \subseteq 2^{V}$ be a 
family of subsets of $V$. A {\em Maker-Breaker game} involves two  
players, named Maker and Breaker, who alternately   
occupy the elements of $V$, called the {\em board} of the game. The Breaker makes the first move. The game
ends when there are no unoccupied elements of $V$. Maker wins the
game if, in the end, the vertices occupied by Maker contain (as a
subset) at least one of the sets in $\mathcal{F}$, the family of {\em
winning sets} of the game. Otherwise Breaker wins.

Chv\'atal and Erd\H{o}s \cite{ChEr} were the first to consider
biased Maker-Breaker games on the edge set of the complete graph.   
They realized that such graph games are often ``easily'' won by
Maker when played fairly (that is when Maker and Breaker each claim
one element at a time). Thus, for many graph games, it is logical to 
give Breaker some advantage. In a $(1:b)$ Maker-Breaker game we
follow the same rules as above, except that Breaker claims $b$ elements each round. It is not too   
difficult to see that, if for some fixed game, Maker can win the $(1:b)$
game, then Maker can win the $(1:b')$ game for every $b'<b$. Therefore Maker-Breaker games are
{\em bias monotone}. It is thus  natural to consider the {\em critical bias} of a game, 
defined as the maximal $b_0$ such that Maker wins the $(1:b_0)$
game.

One of the first biased games that Chv\'atal and Erd\H{o}s
considered in \cite{ChEr} was the Hamiltonicity game played on the
edge set of the complete graph. They proved that Maker wins the $(1:1)$ game, and that for any fixed positive $\varepsilon$ and $b(n)
\ge (1+\varepsilon)\frac{n}{\log n}$, Breaker wins the $(1:b)$ game for large enough $n$. Despite many results by various researchers (see \cite{Krivelevich} for the history), 
the problem of determining the critical bias of this game was open until recently. It was resolved by
Krivelevich \cite{Krivelevich} who proved that the critical bias    
is asymptotically $\frac{n}{\log n}$. We refer the reader
to \cite{Beck4, HKSS} for more information on Maker-Breaker games, as well as
general positional games.

In the spirit of this survey we would like to strengthen Dirac's theorem  from the view point of the Maker-Breaker game.
Let $G$ be a Dirac graph and consider the Hamiltonicity Maker-Breaker game played on $G$. 
Can Maker win the $(1:1)$ game and if yes, what is critical bias?
As the step answering this question, the following theorem \cite{KLS2} established the threshold $b_0$ such that if
$b \ll b_0$, then Maker wins, and if $b \gg b_0$, then Breaker wins.

\begin{theorem} \label{Dirac-games}
There exists a constant $c>0$ such that for $b\le\frac{cn}{\log n}$ 
and a Dirac graph $G$ on $n$ vertices, Maker has a winning strategy for the $(1:b)$
Maker-Breaker Hamiltonicity game on $G$.
\end{theorem}

The theorem implies that the critical bias of this game has order of
magnitude $\frac{n}{\log n}$ (since the critical bias is at most
$(1+o(1))\frac{n}{\log n}$ by the result of Chv\'atal and Erd\H{o}s 
mentioned above). Note that in this theorem, once all the elements  
on the board are claimed, the edge density of Maker's graph is of  
order of magnitude $\frac{\log n}{n}$ and that this is the same as in    
Theorem \ref{Dirac-random-subgraph}. This suggests that as in many other
Maker-Breaker games, the ``probabilistic intuition'', a relation
between the critical bias and the threshold probability of random   
graphs holds here (see \cite{Beck4}). In fact, this
is not a coincidence, and the proofs of Theorems \ref{Dirac-random-subgraph} and \ref{Dirac-games} 
are both done in one  unified framework.

As we already mentioned in Section 2, the random graph 
$G(n,p)$ is a.a.s. Hamiltonian for $p\geq \frac{\log n+\log \log n+\omega(n)}{n}$. Moreover, the
Hamiltonicity threshold coincides with the threshold of having minimum degree $2$.
Can one strengthen these facts using the framework of Maker-Breaker games? Given a graph $G$
with a vertex of degree $3$, Breaker can claim at least two edges incident to this vertex. This will leave 
Maker with a graph containing a vertex of degree $1$ (or less), which is not Hamiltonian. 
Note that even if we allow Maker to move first, it is enough to have two non-adjacent vertices of degree $3$ for Breaker to win.
This will a.a.s. happen in $G(n,p)$ if $p\leq \frac{\log n+3\log \log n-\omega(n)}{n}$. On the other hand, for 
$p\geq \frac{\log n+3\log \log n+\omega(n)}{n}$ the minimum degree of a random graph is a.a.s. $4$ so
it is natural to believe that Maker can now build a Hamilton cycle.
Indeed, such a robust version of the Hamiltonicity of random graphs was proved by Ben-Shimon, Ferber, Hefetz and Krivelevich \cite{BFHK}.

\begin{theorem} \label{random-games}
If $p\geq \frac{\log n+3\log \log n+\omega(n)}{n}$ then a.a.s. Maker has a winning strategy for the $(1:1)$
Hamiltonicity Maker-Breaker game on $G(n,p)$.
\end{theorem}
For larger values of the edge probability $p$ this theorem suggests that the $(1:1)$ game should be an easy win for Maker.
The real question for denser random graphs is therefore to determine the critical bias for the Hamiltonicity game.
Extending the above-mentioned result of Krivelevich,  Ferber, Glebov, Krivelevich and Naor \cite{FGKN} proved that 
for $p \gg \frac{\log n}{n}$ a.a.s. the threshold bias for the Hamiltonicity game on $G(n,p)$
is asymptotically $\frac{np}{\log np}$. An interesting topic for future research is to study Hamiltonicity games 
on random directed graphs.

\section{Compatible Hamilton cycles}
In this section, we consider the following setting which leads to yet another type of robustness that can be used to 
study the Hamiltonicity of Dirac graphs and similar questions. Suppose that we are given a Dirac graph $G$ together with a set of restrictions on its edges given to us by an adversary.
The adversary can forbid certain pairs of edges to appear together on the Hamilton cycle we are trying to build.
Finding a Hamilton cycle satisfying all the imposed restrictions will certainly show some kind of robustness of a
Dirac graph with respect to Hamiltonicity. The type of restrictions 
the adversary can impose is described formally by the following definition. 

\begin{definition}
Let $G=(V,E)$ be a graph.

\vspace{-0.2cm}
\begin{itemize}
  \setlength{\itemsep}{0pt} \setlength{\parskip}{0pt}
  \setlength{\parsep}{0pt}
\item An {\em incompatibility system} $\mathcal{F}$ defined over $G$ is a family $\mathcal{F}=\{F_v\}_{v\in V}$ such that for every $v \in V$, $F_v$ is a family of unordered pairs $F_v
\subseteq \{\{e,e'\}: e\ne e'\in E, e\cap e'=\{v\}\}$.

\item If $\{e,e'\} \in F_v$ for some edges $e,e'$ and a vertex $v$, then we
say that $e$ and $e'$ are \emph{incompatible} in $\mathcal{F}$.
Otherwise, they are \emph{compatible} in $\mathcal{F}$. A subgraph $H \subseteq G$
is \emph{compatible} with $\mathcal{F}$, if all its pairs of edges $e$ and $e'$ are compatible.
\item  For a positive integer $\Delta$, an incompatibility system $\mathcal{F}$
is {\em $\Delta$-bounded} if for each vertex $v \in V$ and an edge $e$ incident to $v$, there are at most $\Delta$ other edges $e'$ incident to $v$ that are
incompatible with $e$.
\end{itemize}
\end{definition}

\noindent
This definition is motivated by two concepts in graph theory.
First, it generalizes \emph{transition systems} introduced by Kotzig \cite{Kotzig68}
in 1968. In our terminology, a transition system is simply a $1$-bounded
incompatibility system. Kotzig's work was motivated by a problem of Nash-Williams on cycle coverings
of Eulerian graphs (see, e.g., Section 8.7 of \cite{Bondy2}).

Incompatibility systems and compatible Hamilton cycles also generalize  
the concept of {\em properly colored} Hamilton cycles in edge-colored graphs.
A cycle is properly colored if its adjacent edges have distinct colors.
The problem of finding properly colored Hamilton cycles in an edge-colored
graph was first introduced by Daykin \cite{Daykin}. He asked whether there
exists a constant $\mu$ such that for large enough $n$, every
edge-coloring of the complete graph $K_{n}$ in which each vertex is incident to
at most $\mu n$ edges of the same color contains a properly colored Hamilton cycle (we refer to
such a coloring as a \emph{$\mu n$-bounded edge coloring}).
Daykin's question has been answered
independently by Bollob\'as and Erd\H{o}s \cite{BoEr76} with $\mu=\frac{1}{69}$,
and by Chen and Daykin \cite{ChDa} with $\mu=\frac{1}{17}$. Bollob\'as
and Erd\H{o}s further conjectured that all $(\lfloor \frac{n}{2}\rfloor-1)$-bounded
edge colorings of $K_{n}$ admit a properly colored Hamilton cycle.
After several subsequent improvements (see, e.g., \cite{AlGu}), Lo \cite{Lo12} recently settled the conjecture asymptotically,
proving that for any positive $\varepsilon$, every $(\frac{1}{2}-\varepsilon)n$-bounded
edge coloring of $K_n$ admits a properly colored Hamilton cycle.

Note that a $\mu n$-bounded edge coloring obviously defines also a $\mu n$-bounded
incompatibility system, and thus the question mentioned above can   
be considered as a special case of the problem of finding compatible
Hamilton cycles. However, in general, the restrictions introduced
by incompatibility systems need not come from an edge-coloring of graphs,
and therefore results on properly colored Hamilton cycles
do not necessarily generalize easily to incompatibility systems.

The study of compatible Hamilton cycles in Dirac graphs, although interesting in its own right, is 
further motivated by the following problem of H\"aggkvist (see \cite[Conjecture 8.40]{Bondy2}).
In 1988 he conjectured that for every $1$-bounded incompatibility
system $\mathcal{F}$ over a Dirac graph $G$, there exists a Hamilton
cycle compatible with $\mathcal{F}$. This conjecture can be settled using Theorem \ref{Dirac-games} on the
Hamiltonicity Maker-Breaker game played on Dirac graphs. Recall that this theorem
asserts the existence of a positive constant $\beta$ such that
Maker has a winning strategy in a $(1:\beta n/\log n)$
Hamiltonicity Maker-Breaker game played on Dirac graphs.
To see how this implies the conjecture, given a graph $G$ and
a $1$-bounded incompatibility system
$\mathcal{F}$, consider a Breaker's strategy
claiming, at each turn, the edges that are incompatible with
the edge that  Maker claimed in the previous turn. This strategy
forces Maker's graph to be compatible with $\mathcal{F}$ at all stages. Since Maker has
a winning strategy for the $(1:2)$ game, we see that there exists a
Hamilton cycle compatible with $\mathcal{F}$. 

Note that the above analysis gives much more, asserting the existence of a compatible Hamilton cycle for every 
$\frac{1}{2}\beta n/\log n$-bounded incompatibility system. Is this the best possible result?
Can we find a compatible Hamilton cycle for every $\Delta$-bounded system, when $\Delta$ is linear in $n$?
The following theorem, proved by Krivelevich, Lee and the author \cite{KLS4}, answers these questions.

\begin{theorem} \label{thm:main}
There exists a constant $\mu>0$ such that the following holds
for large enough $n$. For every $n$-vertex Dirac graph $G$
and a $\mu n$-bounded incompatibility system $\mathcal{F}$ defined
over $G$, there exists a Hamilton cycle in $G$ compatible with $\mathcal{F}$.   
\end{theorem}

\noindent
This shows that Dirac graphs are very robust against
incompatibility systems, i.e., one can find a Hamilton cycle even after
forbidding a quadratic number of pairs of edges incident to each vertex
from being used together in the cycle.
The order of magnitude is best possible since we can
forbid all pairs incident to some vertex from being used together
to disallow a compatible Hamilton cycle.
However, it is not clear what the best possible value of $\mu$      
is and determining it is an interesting open problem.

The proof in \cite{KLS4} provides the existence of a positive constant $\mu$ of approximately
$10^{-16}$ (although no serious attempt was made to optimize it). On the other hand,
the following variant of a
construction of Bollob\'as and Erd\H{o}s \cite{BoEr76}
shows that $\mu$ is at most $\frac{1}{4}$. Let $n$ be an integer of the form
$4k-1$, and let $G$ be an edge-disjoint union of
two $\frac{n+1}{4}$-regular graphs $G_1$ and $G_2$ on the
same $n$-vertex set. Color the edges
of $G_1$ in red and those of $G_2$ in blue. Note that $G$ does not contain
a properly colored Hamilton cycle since any Hamilton cycle of $G$
is of odd length. Let $\mathcal{F}$
be an incompatibility system defined over $G$, where incident edges
of the same color are incompatible. Then there exists a
Hamilton cycle compatible with $\mathcal{F}$ if and only if
there exists a properly colored Hamilton cycle. Since there is no
properly colored Hamilton cycle, we see that there is no Hamilton cycle compatible with $\mathcal{F}$.

It is also natural to use the notion of compatibility 
to strengthen classical results on the Hamiltonicity of random graphs. An obvious question here is given a $\Delta$-bounded incompatibility 
system over $G(n,p)$, when with high probability we can find a compatible Hamilton cycle? 
Since $G(n,p)$ a.a.s. has no Hamilton cycles for $p \ll \frac{\log n}{n}$,
we need to consider values of $p$ above this threshold. Also, since all degrees in our graph are close to $np$, 
it is logical to ask whether we can take $\Delta$ as large as $\mu np$ for some constant $\mu>0$.
The following theorem from \cite{KLS3} summarizes what we know about this problem.

\begin{theorem} \label{thm:sparse_p}
There exists a positive real $\mu$ such
that for $p\gg\frac{\log n}{n}$, the graph $G=G(n,p)$ a.a.s. has
the following property: for every $\mu np$-bounded incompatibility
system defined over $G$, there exists a compatible Hamilton cycle.
Moreover if $p\gg \frac{\log^8 n}{n}$, then we can take 
$\mu=\Big(1-\frac{1}{\sqrt{2}}-o(1)\Big)np$.
\end{theorem}

Recalling the above discussion, this result can be seen as an answer to a generalized version of
Daykin's question. In fact, Theorem \ref{thm:sparse_p} generalizes it in two directions.
Firstly, we replace properly colored Hamilton cycles by compatible
Hamilton cycles and secondly, we replace the complete graph by
random graphs $G(n,p)$ for $p \gg \frac{\log n}{n}$ (note that
for $p=1$, the graph $G(n,1)$ is $K_n$ with probability $1$). 
It is unclear what the best possible value of $\mu$ is. 
The example of Bollob\'as and Erd\H{o}s \cite{BoEr76} of a
$\lfloor \frac{1}{2}n \rfloor$-bounded edge-coloring of $K_n$ with no
properly colored Hamilton cycles implies that the optimal value of
$\mu$ is at most $\frac{1}{2}$, since it provides an upper bound for the case $p=1$.

\vspace{0.2cm}
\noindent
{\bf Sketch of proof of Theorem \ref{thm:sparse_p}.}\, Let $\cal F$ be a $\mu np$-bounded incompatibility
system over $G=G(n,p)$ for some small, but constant, $\mu$. First we construct an expander  graph $R\subset G$ which is
compatible with $\cal F$ and has at most $O(n)$ edges and such that every $X\subset V(R)$ of size at most $n/4$ has at least $2|X|$ neighbors outside $X$. This graph can be obtained by choosing, randomly with replacement, $d$ (a large constant) edges incident to every vertex in $G$ and showing that the resulting graph has all the desired properties with positive probability. Once we have $R$, consider a path $P$ in $G$ such that $P \cup R$ is compatible with $\cal F$ and $P$ is a longest path in $P \cup R$. Using the rotation-extension technique (from the previous section) together with the properties of $R$, we can rotate both endpoints of $P$  to obtain $\Omega(n^2)$  pairs $(u,v)$ such that we can extend $P$ by adding the edge $(u,v)$. Let us call these pairs {\em boosters}. Since $G=G(n,p)$ is a random graph, a.a.s. $\Omega(n^2p)$ of the boosters are actual edges in $G$.
Moreover, since $R$ has only $O(n)$ edges, one can show that one of the boosters is compatible with $P\cup R$. We can therefore repeat this process until we obtain a compatible Hamilton cycle.
\hfill $\Box$
\vspace{0.2cm}

To complete this section, we mention another closely related way to get a robust version of the 
Hamiltonicity result for random graphs. In an edge-colored graph, we say that a subgraph is {\em rainbow} if all its edges have distinct colors. There is a vast literature on the branch of Ramsey theory where one seeks rainbow
subgraphs in edge-colored graphs. Note that one can easily avoid rainbow copies by using the same color for all edges, and hence in order to find a rainbow subgraph one usually imposes
some restrictions on the distribution of colors. Erd\H{o}s, Simonovits and S\'os \cite{ErSiSo} and Rado \cite{Rado} developed anti-Ramsey theory where one attempts to
determine the maximum number of colors that can be used to color the edges of the complete graph without creating a rainbow copy of a fixed graph. In a different direction, one can try to
find a rainbow copy of a target graph by imposing global conditions on the coloring of the host graph. For a real $\Delta$, we say that an edge-coloring of $G$ is {\em globally
$\Delta$-bounded} if each color appears at most $\Delta$ times on the edges of $G$. In 1982, Erd\H{o}s, Ne\v set\v ril and R\"odl \cite{ENR} initiated the study of the problem of finding
rainbow subgraphs in a globally $\Delta$-bounded coloring of graphs. One very interesting question of this type is to find sufficient conditions for the existence of a rainbow Hamilton cycle
in any globally $\Delta$-bounded coloring. Substantially improving earlier results, Albert, Frieze and Reed \cite{AlFrRe} proved the existence of a
constant $\mu>0$ for which every globally $\mu n$-bounded coloring of $K_n$ (for large enough $n$) admits a rainbow Hamilton cycle. In fact, they proved a stronger statement, asserting
that for all graphs $\Gamma$ with vertex set $E(K_n)$ (the edge set of the complete graph) and maximum degree at most $\mu n$, there exists a Hamilton cycle in $K_n$ which is also an
independent set in $\Gamma$.

It turns out that the proof technique used in proving Theorem \ref{thm:sparse_p} can be easily modified to give the following result (see \cite{KLS3}) that extends the above to random graphs.

\begin{theorem} \label{thm:subramsey}
There exists a constant $\mu > 0$ such that for $p \gg \frac{\log n}{n}$,
the random graph $G=G(n,p)$ a.a.s.~has the following property: 
every globally $\mu np$-bounded coloring of $G$ contains a rainbow Hamilton cycle.
\end{theorem}

Theorem \ref{thm:subramsey} is clearly best possible up to the constant $\mu$ since one can
forbid all rainbow Hamilton cycles in a globally $(1+o(1))np$-bounded coloring by simply
coloring all edges incident to some fixed vertex with the same color.

\section{Resilience of graph properties}
In this section we discuss the following question.
Given a graph $G$ with some property $\cal P$, how much
should one change $G$ in order to destroy $\cal P$? 
We call this {\it the resilience of $G$ with respect to $\cal P$}.
There are two natural kinds of resilience: global and local. It is
more convenient to first define these quantities with respect to
monotone increasing properties. Recall that $\cal P$ is monotone increasing if it
is preserved under edge addition.

\begin{definition}
Let $G$ be a graph and $\cal P$ be a monotone increasing property. 
\begin{itemize}
\item 
The global resilience  of $G$ with respect to $\cal P$ is the minimum number $r$
such that one can obtain a graph not having
$\cal P$ by deleting $r$ edges from $G$. 
\item
The local resilience of a graph $G$ with respect to $\cal P$ is the minimum number $r$ such
that one can obtain a graph not having $\cal P$ by deleting at most $r$ edges at every vertex of $G$. 
\end{itemize}
\end{definition}

\noindent 
The notion of global resilience  is not new. In fact,
problems about global resilience  are  popular in extremal graph
theory. For example, the celebrated Tur\'an problem can be rephrased as follows.
How many edges should one delete from the complete graph $K_n$ to make it $H$-free for some fixed graph $H$?
The notion of local resilience is more recent and its systematic study was initiated by the author together with Vu \cite{SV}.
It is motivated by an observation that many interesting properties are easy to destroy by small local changes.
For example, to destroy the Hamiltonicity  it suffices to delete all edges incident to one vertex, which is much less than  
the total number of edges in a graph.

If a graph property is not monotone, like containing an induced copy of a fixed graph or having a
trivial automorphism group, then we may have to
both delete and add edges. This leads  to the following more general definition.
\begin{definition}
Given a property $\cal P$. The global/local resilience  of $G$ with respect
to $\cal P$ is the minimum number $r$ such that there is a graph $H$ on
$V(G)$ with the total number of edges/maximum degree at most $r$ for which the graph $G
\triangle H$ does not have $\cal P$. 
\end{definition}

One can observe that there is a certain duality between properties
and resilience. If the property is local (such as containing a
triangle), then it makes more sense to talk about the global
resilience. On the other hand, if the  property is global (such as being
Hamiltonian), then the local resilience  seems to be the right
parameter to consider. There are some interesting exceptions to this rule, like chromatic 
number, which we discuss later in this section.

Using the above notions one can easily generate many questions by choosing some
graph with an interesting property and asking for its resilience with respect to this property. 
For example, asking for the local resilience of the complete graph with
respect to Hamiltonicity leads immediately to the celebrated theorem of Dirac, which 
has been discussed in depth earlier in this paper. Following this approach, in the rest of 
this section we will revisit several classical theorems on random graphs 
and discuss the corresponding resilience results.

\subsection{Perfect Matching and Hamiltonicity}
Let $G$ be a graph on $n$ vertices, where $n$ is even.
A perfect matching in $G$ is a set of disjoint edges that covers all $n$
vertices. The threshold for the appearance of such a matching in $G(n,p)$ was determined by Erd\H{o}s and R\'enyi  
already in one of their early papers on random graphs. In \cite{ER2} they proved
that if $p \geq  \frac{\log n+ \omega(n)}{n}$, then a.a.s. $G(n,p)$ has a perfect matching. 
This determines the threshold since for $p \leq  \frac{\log n- \omega(n)}{n}$ we a.a.s. have isolated vertices.

To give a lower bound for the local resilience  of 
having a perfect matching consider the following way to destroy
this property. Split the vertex set of
$G$ into two parts $X$ and $Y$ of size $n/2+1$ and $n/2-1$
respectively. Then delete all edges inside the set $X$. Thus $X$
becomes an independent set  and it is impossible to match all of
its vertices  with vertices of $Y$, since $|Y|<|X|$. 
In $G(n,p)$ (with $p$ sufficiently large), with high probability
all vertices have degree $(1+o(1))np$. Thus,  one would expect
that a.a.s. a random graph has an induced subgraph on $n/2+1$ vertices whose maximum degree is $(1/2+o(1))np$. 
So, the local resilience  of $G(n,p)$ with respect to having a perfect matching is a.a.s. at most
$(1/2+o(1))np$.

It was proved in \cite{SV} that for $p \gg \log n/n$ this trivial upper
bound is actually the truth. Moreover, in this case one can obtain the following rather accurate estimates for error terms.
If $G'$ is a subgraph of $G=G(n,p)$ with maximum degree at most $np/2-8\sqrt{np\log n}$ then a.a.s. $G-G'$ has a perfect matching. On the other hand, $G(n,p)$ a.a.s. contains a subgraph $G''$ with maximum degree at most $np/2+2\sqrt{np\log n}$
such that $G-G''$ has no perfect matching. The proof of the first part is not very difficult. By randomly splitting the vertices of $G-G'$ into two sets of size $n/2$ we can obtain a bipartite subgraph $H\subset G-G'$ with very good expansion properties. Using these properties one can verify that $H$ satisfies Hall's condition (see \cite{Diestel}), i.e., every subset $X$ on one side has at least $|X|$ neighbors on the other side. The second part follows by analyzing more carefully the degree distribution of the subgraph of $G(n,p)$ on $n/2+1$ vertices.

The question of resilience of random graphs with respect to Hamiltonicity is substantially  more difficult.
Here too one can destroy all Hamilton cycles of $G=G(n,p)$ by splitting the vertex set of $G$ into
two parts whose sizes differ by at most two and deleting all the
edges inside the larger part. Therefore, similar as above, the local resilience  of the random graph with respect to Hamiltonicity
is a.a.s. at most $(1/2+o(1))np$. Note that we can also destroy Hamiltonicity by simply disconnecting the  graph.
To do so,  split the vertex set into two parts whose sizes differ by at most one  and delete all edges between
them. However in the case of random graphs, this does not change the asymptotics of our lower
bound for resilience, since bipartite graph we remove has typically maximum degree $(1/2+o(1))np$.
This suggests that the local resilience of $G(n,p)$ with respect to
Hamiltonicity is with high probability $(1/2+o(1))np$, at least when the random graph is reasonably dense.
This was indeed proved for $p \geq \log^4 n/n$ in \cite{SV}, where it was conjectured that 
the value of edge probability $p$ can be decreased all the way to $p \gg \log n/n$.
This conjecture was resolved by Lee and Sudakov \cite{LS}, who proved the following result.

\begin{theorem} \label{theorem:Hal1} 
For every $\varepsilon>0$, there exist a constant $C=C(\varepsilon)$ such that if
$p \geq C\log n/n$, then a.a.s. every subgraph $G$ of $G(n,p)$ with minimum degree at least $(1/2+\varepsilon)np$ is Hamiltonian.
\end{theorem}
\noindent
Since a complete graph on $n$ vertices is also a random graph $G(n,p)$ with $p=1$,
this theorem can be viewed as a far reaching generalization of Dirac's theorem.

The proof of Theorem \ref{theorem:Hal1} uses the rotation-extension technique, mentioned in Section 4.1, with a few additional ideas. One is to split the graph $G$ into two graphs $G_1, G_2$, where the first graph has only small fraction of the edges and will be used to perform rotations and the second graph $G_2$ will be used to perform extensions. We also partition the longest path $P$ in $G_1$ into several intervals $I_1, \ldots, I_k$ and show that for most indices $j$ the majority of newly constructed paths have their broken edges outside the interval $I_j$. Therefore these paths traverse $I_j$ either in the original or the reverse order. This is used to show that, rotating the longest path in $G_1$, one can obtain a set $S_P$ of $(1/2+\varepsilon)n$ new endpoints. Furthermore for every $v \in S_P$, rotating the path again and keeping $v$ fixed we can obtain a set $T_v$ of size at least $(1/2+\varepsilon)n$ such that for all $w \in T_v$ there is a longest path in $G_1$ starting at $v$ and ending in $w$. Since $G_2$ contains most of the edges of $G$ and the minimum degree of $G$ is at least $(1/2+\varepsilon)np$ one can show that one of these paths can be closed into cycle by using the edges of $G_2$. Repeating this procedure several times we obtain a Hamilton cycle.

When the edge probability $p$ is close to the threshold for Hamiltonicity, i.e., $p=(1+\varepsilon)\log n/n$, the random graph $G(n,p)$  has some vertices whose degrees are 
significantly smaller than $np$. Therefore in this range of edge probability, one can easily create isolated vertices by deleting  $(1/2+o(1))np$ edges incident to such vertices. 
Hence in this case, we need to revise the definition of resilience and make the number of edges we allow to delete at a vertex $v$ depending on the degree of $v$. For further discussion of this regime and known results see e.g., \cite{BKS} and its references.

Many techniques developed for Hamiltonicity of random graphs 
rely only on properties of the edge distribution of $G(n,p)$ and therefore 
can be used to study pseudo-random graphs as well. There are several closely related definitions
of pseudo-random graphs (see, e.g., \cite{KSu} for discussion). 
Here we use the following one, which is based on spectral properties of such graphs.
Consider a graph $G$ on $n$ vertices. Since its adjacency matrix is symmetric, it has $n$ real eigenvalues which we denote by   $\lambda_1(G) \geq \lambda_2 (G) \geq
\ldots \geq \lambda_n (G)$. The quantity $\lambda(G)= \max_{i \geq 2} |\lambda_i(G)|$
is called the {\em second eigenvalue} of $G$ and plays an
important role. We say that $G$  is  an {\em
$(n,d,\lambda)$-graph} if it is $d$-regular, has $n$ vertices and
$\lambda(G)$ is at most $\lambda$. 

It is well known (see, e.g.,
\cite{KSu}) that if $\lambda$ is much
smaller than the degree $d$, then $G$ has strong pseudo-random
properties, i.e., the edges of $G$ are distributed like in the random
graph $G(n,d/n)$. Therefore it is natural to ask whether  (similar to $G(n,d/n)$)
pseudo-random graphs with $\lambda \ll d$ are Hamiltonian. 
Such a result  was establish by the author together with Krivelevich in \cite{KSu2}, where 
we proved that $(n,d,\lambda)$-graphs with $\lambda<d/\log n$ (see \cite{KSu2} for slightly better bound)
are Hamiltonian. It would be very interesting to determine the right order of magnitude of 
the ratio $d/\lambda$ which already implies Hamiltonicity. Together with Krivelevich \cite{KSu}
we proposed the following conjecture.
\begin{conjecture}
	There exists a positive constant $C$ such that for all sufficiently large $n$, every 
	 $(n,d,\lambda)$-graph with $d/\lambda>C$ contains a Hamilton cycle.
\end{conjecture}
One can further show that when $d \gg \lambda$ the corresponding graph is robustly Hamiltonian. Indeed, together with Vu \cite{SV} we proved that if $d/\lambda>\log^2 n$ (actually $\log^{1+\delta} n $ is enough) then the local resilience of such an $(n,d,\lambda)$-graph with respect to Hamiltonicity is
$(1/2+o(1))d$. By the above discussion the constant $1/2$ is best possible. 

Another interesting open problem is the resilience of random regular graphs with respect to Hamiltonicity.
The $n$-vertex {\em random $d$-regular graph}, which is denoted by $G_{n,d}$ (where $dn$ is even), is the uniform probability space of all
$d$-regular graphs on $n$ vertices labeled by the set $[n]$. In this model, one cannot apply the
techniques used to study $G(n, p)$ as these two models do not share the same probabilistic
properties. Whereas the appearances of edges in $G(n, p)$ are independent, the appearances
of edges in $G_{n,d}$ are not. Nevertheless, many results obtained thus far for the random
regular graph model $G_{n,d}$ are very similar to the results obtained in $G(n, p)$
with suitable expected degrees, namely, $d = np$. For a detailed discussion of random regular graphs and their properties we refer the interested reader to the
excellent survey of Wormald \cite{Wor}. 

It was proved by Robinson and Wormald \cite{RobWor} that the random regular graph $G_{n,d}$ is a.a.s. Hamiltonian already for 
$d \geq 3$. Therefore the question of local resilience of $G_{n,d}$ with respect to this property already makes sense for constant degrees.
By the above discussion, it is logical to guess that this resilience should be typically of order $d/2$, at least for large $d$. 
More precisely, together with Ben-Shimon and Krivelevich \cite{BKS2} we made the following conjecture.
\begin{conjecture}
	For every $\varepsilon > 0$ there exists an integer $d_0=d_0(\varepsilon)$ such that, for every fixed
	integer $d > d_0$, asymptotically almost surely the local resilience of $G_{n,d}$ with respect to Hamiltonicity is at least 
	$(1/2-\varepsilon)d$ and at most $(1/2+\varepsilon)d$.
\end{conjecture}

\noindent
The upper bound of this conjecture follows easily from the known properties of the edge distribution of $G_{n,d}$ together with constructions mentioned above.
So the main task is to prove the lower bound. When $d\gg \log^2 n$ one can prove this conjecture using the above mentioned result from \cite{SV} on $(n,d,\lambda)$-graphs
and the known properties of the second eigenvalue of $G_{n,d}$ (see, e.g., \cite{Wor, KSVW}). For smaller values of $d$ this conjecture is still open, although one can show (see \cite{BKS2}) that a.a.s. the local resilience is linear
in $d$.

\subsection{Chromatic number}
One of the most important parameters of the 
random graph $G_{n,p}$ is its
chromatic number, which we denote by 
$\chi(G_{n,p})$. Trivially for every graph
$\chi(G)\geq |V(G)|/\alpha(G)$, where $\alpha(G)$
denotes the size of a largest independent set in $G$. 
It can be easily shown, using first moment computations,
that a.a.s. $\alpha(G_{n,p})\leq 2\log_b (np)$, where
$b=1/(1-p)$ (all other logarithms in this paper are in the natural base $e$). 
This provides a lower bound on the chromatic number of the random graph, 
showing that with high probability $\chi(G_{n,p})\geq \frac{n}{2\log_b (np)}$. 
The problem of determining the asymptotic behavior  of $\chi(G_{n,p})$, 
posed by Erd\H{o}s and R\'enyi in the early 60s, stayed for many years 
as one
of the major open questions in the theory of random graphs until
its solution by Bollob\'as \cite{Bol2}, using a
novel approach based on martingales that enabled him to prove that
a.a.s. $\chi(G_{n,p})=(1+o(1)) \frac{n}{2\log_b (np)}$  for dense random graphs. Later \L uczak \cite{L} showed that this estimate also holds for all values of $p \geq c/n$. 

Given two arbitrary graphs $G$ and $H$ it is a folklore result which is easy to prove (see,
e.g., \cite{Lovaszbook} Chapter 9) that
$$\chi (G \cup H) \le \chi (G) \chi (H).$$
Moreover there are pairs of graphs for which the equality  holds.
Therefore adding a few edges or low degree graph to a graph $G$ sometimes can have  
a substantial impact on its chromatic number. The question, whether this is the case
for random graphs was first posed in \cite{SV}, where the authors study the resilience of the chromatic number of $G(n,p)$. Unlike for Hamiltonicity, changing the adjacency of few vertices in a graph usually has only a minor effect on its chromatic number. Therefore it is interesting to study both global and local resilience of the chromatic number.

Note that for an arbitrary graph $G$ one can double its chromatic number by choosing an arbitrary subset of 
$2\chi(G)$ vertices and adding all the missing edges to make this subset a clique.
This gives an upper bound $O(\chi^2(G))$ on the global resilience and a bound $O(\chi(G))$ on a local resilience of the chromatic number. 
Together with Vu \cite{SV} we conjectured that this is with high probability tight for relatively dense random graphs (say $p \geq n^{-1+\delta}$ for any fixed $\delta>0$).  Improving the earlier results from
\cite{SV}, Alon and Sudakov \cite{AS} proved the following theorem, which makes a substantial progress on this conjecture. 

\begin{theorem}
	\label{resilience-chromatic}
	Let $\varepsilon>0$ be a fixed constant and let $n^{-1/3+\delta} \leq p \leq 1/2$ for some $\delta>0$. Then a.a.s.
	\begin{enumerate}
		\item
		for every collection $E$ of $2^{-12}\varepsilon^2\frac{n^2}{\log^2_b (np)}$ edges the chromatic number of $G_{n,p} \cup E$ is still 
		at most $(1+\varepsilon)\frac{n}{2\log_b (np)}$,
		\item
		and for every graph $H$ on $n$ vertices with maximum degree 
		$\Delta(H) \leq 2^{-8}\varepsilon \frac{n}{\log_b (np) \log \log n}$ the chromatic number of $G_{n,p} \cup H$ is still
		at most $(1+\varepsilon)\frac{n}{2\log_b (np)}$.
	\end{enumerate}
\end{theorem}
This determines the global resilience of dense random graphs up to a constant factor and leaves only a small multiplicative gap of $\log \log n$ for the local resilience. Both these results show that adding quite large and dense graphs to $G_{n,p}$ with $ p \gg n^{-1/3}$ typically has very little impact on its chromatic number.
For $p$ below $n^{-1/3}$ much less is known. It was proved in \cite{SV} that 
for every positive integer $d$ and for every $\varepsilon >0$ there is a constant
$c=c(d,\varepsilon) $ such that the following holds. For all $p > c/n$, 
adding to $G(n,p)$ any graph with maximum degree $d$ a.a.s. cannot increase its chromatic number by a factor 
of larger than $(1+\varepsilon)$.

\vspace{0.2cm}
\noindent
{\bf Sketch of proof of Theorem \ref{resilience-chromatic} (1).}\, Choose $k_0=(2-o(1)) \log^2_b (np)$ such that the expected number 
of independent sets of this size in $G(n,p)$ is $\mu={n \choose k_0}(1-p)^{k_0 \choose 2}>n^4$. The expected number of such independent sets containing a given pair of vertices is  
$$\mu_0={n-2 \choose k_0-2}(1-p)^{k_0 \choose 2}=(1+o(1))\frac{k_0^2}{n^2}.$$ Let $\cal I$ be a largest collection of independent sets of size 
$k_0$ in $G(n,p)$ such that no pair of vertices belongs to more than $4\mu_0$ of these sets. This condition shows that by changing one edge of the random graph we cannot affect the size of $\cal I$ by more than $4\mu_0$. Therefore using a standard estimate for the tails of martingales we can show that $|\cal I|$ equals $(1-o(1))\mathbb{E}[|\cal I|]$
with probability at least $1-2^{-2n}$. Moreover one can also show (see \cite{KSVW2}) that $\mathbb{E}[|{\cal I}|]=(1-o(1))\mu$.

Let  $E$ be the collection of  $2^{-12}\varepsilon^2\frac{n^2}{\log^2_b (np)}$ edges which was added to $G(n,p)$. Consider an auxiliary bipartite graph $H$ with parts $\cal I$ and $E$ in which an independent set $I \in {\cal I}$ is adjacent to an edge $(u,v) \in E$ iff
both vertices $u,v$ belong to $I$. By the definition of $\cal I$, every edge $(u,v) \in E$ is contained in at most $4\mu_0$ sets from $\cal I$.
Therefore the number of edges $e(H)$ is bounded by $4\mu_0 m$. Thus there is an independent set $I \in {\cal I}$, whose degree in $H$ is at most
$e(H)/|{\cal I}|$. Such $I$ contains at most
$$\frac{e(H)}{|{\cal I}|}\leq \frac{4\mu_0 |E|}{|{\cal I}|}\leq \frac{5\mu_0 |E|}{\mu}\leq 6k_0^2\frac{|E|}{n^2} \leq 2^{-7}\varepsilon^2$$ 
edges from $E$. 

This argument shows that we can find an independent set in $G(n,p)$ which contains no edges from $E$. Color it by color one, remove it from the graph and continue this process.
One can further prove that, as long as the number of remaining vertices is larger than $\varepsilon n/\log_b (np)$, we can still find an independent set $I$ 
in $G(n,p)$ which has very few edges from $E$. Then (by a standard lower bound on the independence number in a graph with a given number of edges)  there is a large independent set in $I$ not containing any edges of $E$ which we can color by a new color. Using some simple computations, see \cite{AS}, we can then show that the total number of colors used in this process is at most $(1+O(\varepsilon))\frac{n}{2\log_b (np)}$. Coloring the  
remaining $\varepsilon n/\log^2_b (np)$ vertices by additional colors we obtain a coloring of $G(n,p) \cup E$ into 
$(1+O(\varepsilon))\chi(G(n,p))$ colors.
\hfill $\Box$

\subsection{Symmetry}
An {\em automorphism} of a graph $G$ is a permutation $\pi: V(G) \rightarrow V(G)$ of the vertices of $G$
such that $\big(\pi(u),\pi(v)\big)$ is an edge of $G$ if and only if $(u,v)$ is an edge of $G$. The collection of all automorphisms of $G$ forms a group which is denoted by $Aut(G)$. It is clear that for any
graph $G$ the identity belongs to $Aut(G)$. We say that $Aut(G)$ is
trivial or equivalently $G$ is {\em asymmetric} if $Aut (G)$ does not contain any
permutation other than the identity, otherwise we call $G$ {\em symmetric}.
The automorphism group was one of the first objects studied by P. Erd\H os and A. R\'enyi in their sequence of papers which started the theory of random graphs. In 1963 they proved  that for $p \ge (1+\ep) \log n /n $,  $G(n,p)$ is asymptotically almost surely asymmetric. In fact, Erd\H os and R\'enyi studied
a more general question of how much should one change the random graph to have a non-trivial
automorphism. In \cite{ER1} they proved that if both $1-p, p \gg \log n/n$, then a.a.s. we need to add and delete at least
$(2+o(1))np(1-p)$ edges from $G(n,p)$ to obtain a symmetric graph. This determines 
the global resilience of $G(n,p)$ with respect to having a trivial automorphism group.

To explain the quantity $(2+o(1))np(1-p) $, notice that next to the identity, the simplest permutation is the
transposition of two vertices, say $u$ and $v$. This permutation is an automorphism if $u$ and $v$ have exactly the same neighbors. To achieve this one needs to delete $d(u) + d(v) -2codeg (u,v)$ edges of $G$, where 
$codeg (u,v)$ is the number of common neighbors of these vertices. It shows that
in the case of $G(n,p)$, with high probability it is enough to delete $2(1+o(1))np(1-p)$
edges to have a symmetric graph. Together with
Kim and Vu, the author extended the above theorem of Erd\H os and R\'enyi  
showing that small local changes are not enough to make a random graph symmetric. In \cite{KSV} they proved the 
following result.

\begin{theorem} If both $1-p, p \gg \log n/n$, then the local resilience
of $G(n,p)$ with respect to being asymmetric
is a.a.s. $(1+o(1))np(1-p) $.
\end{theorem}

It is worth mentioning that the local resilience was not the main object of study in \cite{KSV}, but
it was used as a tool  to prove a new result. As a corollary of the above theorem (more precisely of its proof), it was shown there that a random regular graph of relatively large degree is asymptotically almost surely non-symmetric, confirming a conjecture of Wormald \cite{Wor}. The main idea in
\cite{KSV} was as follows. Consider the indicator graph function
$I(G)$, where $I(G)=1$ if $G$ is non-symmetric and $0 $ otherwise. 
We want to show that with high probability $I(G)=1$ where
$G$ is a random regular graph. One may want to view this
statement as a sharp concentration result, namely, $I(G)$ is a.a.s.
close to its mean. However, it is impossible to
prove a sharp concentration result for a  random variable having
only two values close to each other. The idea here is to ``blow up"
$I(G)$ using the notion of local resilience. Instead of $I(G)$ we
used a function $D(G)$ which (roughly speaking) equals the local
resilience  of $G$ with respect to being non-symmetric. This
function is zero if $G$ is symmetric and rather large otherwise.
This gives us room to show that $D(G)$ is strongly concentrated around a large
positive value, and from this we can conclude that asymptotically almost surely
the random regular graph is non-symmetric.

\subsection{Further results}
There are some additional results on resilience of random graphs whose details we will not discuss in this survey. For convenience of the reader we conclude this section with a short list of such results together with relevant references.

\begin{itemize}
\item
In  the  past few  years,  there  has  been  some  considerable  success  in  extending  classical
results in extremal combinatorics to sparse random settings due to the breakthroughs made independently by Conlon and Gowers \cite{CG} and Schacht \cite{Sch}.
In particular, they proved Tur\'an-type theorems for random graphs, establishing the global resilience with respect to containing a fixed non-bipartite graph $H$.

\item 
Resilience of $G(n,p)$ with respect to {\em pancyclicity}, which is a property of containing cycles of all lengths from $3$ to $n$, was established in \cite{KLS5}.
\item
In \cite{BCS} the authors studied the resilience of the property of
containing given almost spanning tree of bounded degree.
\item
Resilience of a random directed graph with respect to Hamiltonicity was studied in \cite{HSS, FNNPS}.
\item
In \cite{DKMS} the authors studied the resilience of $G(n,p)$ with respect to containing cycle of a given linear length.
\item
An $H$-factor in a graph $G$ is a collection of vertex disjoint copies of a fixed graph $H$ covering all vertices of $G$.
The resilience with respect to containing such a factor was studied in \cite{HLS, BLS}. More generally, 
the resilience with respect to containing almost spanning and 
spanning bounded degree graphs was studied in \cite{BKT, NS, ABET}.
\item
In \cite{FKN} the authors made a very interesting connection between the local resilience in random graphs and winning strategies in biased Maker-Breaker games. 
\end{itemize}

\section{Conclusion}
Although we have made an effort to provide a systematic coverage of
recent robust extensions of various classical results in graph theory and random graphs, there are certainly quite a few of them that were left out of this survey, due to the limitations of space
and time (and of the author's energy). Still, we would like to believe that we have presented enough examples demonstrating how one can revisit known results and use various measures of robustness to ask new and interesting questions. Therefore we hope that this survey will motivate future research on the subject. 

\vspace{0.25cm}
\noindent
{\bf Acknowledgement.} \, We would like to thank F. Dr\"axler, A. Ferber, N. Kam\v{c}ev, M. Krivelevich, M. Kwan, J. Noel, A. Pokrovskiy and the referee for helpful comments.

\end{document}